\documentclass[11pt]{article}
\usepackage{amssymb}
\usepackage{amsmath}
\usepackage{amsthm}
\usepackage{latexsym}
\usepackage{amsfonts}
\usepackage{graphicx}
\usepackage{graphics}

\newcommand{\lyxaddress}[1]{
\par {\raggedright #1
\vspace{1.4em}
\noindent\par}
}

\newcommand{\dis}{\displaystyle}
\theoremstyle{plain}
\newtheorem{thm}{Theorem}[section]   % Αρίθμηση συνεχόμενη (όχι κατά Θεώρημα, Λήμμα κ.λπ.)
\newtheorem{prop}[thm]{Proposition}
\newtheorem{lem}[thm]{Lemma}

\theoremstyle{definition}
\newtheorem{rem}[thm]{Remark}

\newtheorem*{Proof}{Proof}

\newcommand{\lra}{\;\longrightarrow\;}
\newcommand{\ra}{\;\rightarrow\;}

\newcommand{\Ga} {{\varGamma}}

\newcommand{\de}{\delta }
\newcommand{\OO} {{\varOmega}}

\newcommand{\e}{\varepsilon }

\newcommand{\z}{\zeta }

\newcommand{\la}{\lambda }
\newcommand{\mi}{\mu }

\newcommand{\ti}{\tau }

\newcommand{\C}{\mathbb{C}}

\newcommand{\N}{\mathbb{N}}

\newcommand{\oO}{\overline{\varOmega}}

\newcommand{\tde}{\widetilde{\delta}}

\newcommand{\tp}{\widetilde{p}}

\newcommand{\tg}{\widetilde{g}}

\newcommand{\tq}{\widetilde{q}}

\newcommand{\cf}{{\cal{F}}}

\newcommand{\cD}{{\cal{D}}}

\newcommand{\cb}{{\cal{B}}}

\newcommand{\cu}{{\cal{U}}}

\newcommand{\ld}{\ldots}

\newcommand{\sm}{\smallsetminus}

\newcommand{\qb}{$\quad\blacksquare$}
\begin{document}
\title{\bf Generic approximation of functions\\ by their Pad\'{e} approximants}
\author{G. Fournodavlos$^1$ and V. Nestoridis$^2$}

\maketitle

\lyxaddress{Department of Mathematics, University of Athens, Panepistemioupolis,
15784, Athens, Greece, e-mail: $^1$gregdavlos@hotmail.com, $^2$vnestor@math.uoa.gr}
%
%\begin{abstract}
\noindent
%
%{\bf Abstract}.\medskip
\begin{abstract}

Generic approximation of entire functions by their Pad\'{e} approximants has been
achieved in the past (\cite{3}). In the present article we obtain generic approximation of holomorphic functions
on arbitrary open sets by sequences of their Pad\'{e} approximants. Similar results hold with functions
smooth on the boundary of their domain of definition. In addition, the approximation is valid simultaneously
with respect to all centers of expansion.
\end{abstract}
{\em Subject Classification MSC2010}\,: primary 41A21, 30K05
secondary 30B10, 30E10, 30K99, 41A10, 41A20. \vspace*{0.2cm} 
\section{Introduction}\label{sec1} % 1
\noindent

Every holomorphic function on a disc can be approximated by the partial sums of its Taylor development.
Furthermore, on any simply connected domain generically all holomorphic functions are
limits of some subsequence of their partial sums, in the topology of uniform
convergence on compacta (\cite{12}).

The partial sums of Taylor expansions are trivial cases of Pad\'{e} approximants, which are rational functions in general.
In \cite{3} it is proven that generically every entire function $f$ is the limit of a sequence
of its Pad\'{e} approximants, $[p_n/q_n]_f$ , provided that $p_n\ra+\infty$ and $p_n-q_n\ra+\infty$. This answered
positively Baker's conjecture (\cite{1}) for quasi all entire functions.

In the present article we strengthen this result in the following ways. Firstly, we show that alone the condition
$p_n\ra+\infty$ is sufficient. Secondly, we extend this result in any simply connected open set $\OO\subseteq\C$.
Finally, we obtain simultaneous approximations by the Pad\'{e} approximants $[p_n/q_n]_{f,\z}$ with respect to
any center of expansion $\z\in\OO$; see also \cite{12}, \cite{15}. Next, we achieve the same results with functions smooth on the boundary of $\OO$, $f\in A^\infty(\OO)$, provided $(\C\cup\infty)\sm\oO$ is connected (\cite{11}, \cite{13}).

When we make approximations by using polynomials, the maximum principle leads to compact sets with connected
complement. However, Pad\'{e} approximants have in general poles. This allows us to approximate also on compact sets
of arbitrary connectivity. Thus, we obtain generic approximation of holomorphic functions $f\in H(\OO)$, defined on arbitrary open sets, by their Pad\'{e} approximants $[p_n/q_n]_f$. The condition now is $p_n\ra+\infty$ and $q_n\ra+\infty$. Again the approximation is simultaneous with respect to all centers. Lastly, we transfer the same results in the space $A^\infty(\OO)$, as in the previous case.

Another generic approximation which is well known, is universal approximation outside the domain of definition.
Universal Pad\'{e} approximants have been introduced in \cite{5} and although our methods were inspired from the methods used in \cite{5}, they can be applied to universal Pad\'{e} approximants in order to give stronger and more complete results.
This will be done hopefully in the future, in another article.

Our method of proof is based on Baire's Category theorem. Despite the fact that the results may not appear so
elementary or simple, the proofs are. For the role of Baire's theorem in Analysis we refer to \cite{9} and \cite{10}.
Parts of the present article can be found in \cite{6}, \cite{7}.

\section{Preliminaries}\label{sec2}
\noindent

Let $f(z)=\dis\sum^\infty_{v=0}a_v(z-\z)^v$, $a_v\in\C$, be a formal series, $\z\in\C$.
A Pad\'{e} approximant $[p/q]_{f,\z}$ of $f$, $p,q\in\{0,1,2,\ld\}$, is a rational
function of the form
\[
\frac{\dis\sum^p_{v=0}n_v(z-\z)^v}{\dis\sum^q_{v=0}d_v(z-\z)^v}, \ \ d_0=1.
\]
such that its Taylor series $\dis\sum^\infty_{v=0}b_v(z-\z)^v$ coincides
with $\dis\sum^\infty_{v=0}a_v(z-\z)^v$ up to the first $p+q+1$ terms;
that is $b_v=a_v$ for $v=0,\ld,p+q$ (\cite{2}). Thus, if $f$ is a holomorphic
function in a neighborhood of $\z$, we can define $[p/q]_{f,\z}$ corresponding
to the Taylor development of $f$ around $\z\in\C$.

We notice that in case $q=0$ there exists always a unique
Pad\'{e} approximant of $f$ and $[p/q]_{f,\z}(z)= S_p(z)$, where $S_p(z)= \dis\sum^p_{v=0}a_v(z-\z)^v$.
For $q\ge 1$ such rational function does not always exist, neither it is unique, if it exists.
However, it is true that there exists a unique Pad\'{e} approximant of $f$, $[p/q]_{f,\z}$,
if and only if the following Hankel determinant is not zero:
\[
\det\left|\begin{array}{cccc}
  a_{p-q+1} & a_{p-q+2} & \cdots & a_p \\
  a_{p-q+2} & a_{p-q+3} & \cdots & a_{p+1} \\
  \vdots & \vdots &  & \vdots \\
  a_p & a_{p+1} & \cdots & a_{p+q-1} \\
\end{array}\right|\neq0, \ \ a_i=0, \ \ \text{when} \ \ i<0.
\eqno(\ast)
\]
Then we write $f\in \cD_{p,q}(\z)$. Under this condition $[p/q]_{f,\z}$
is given by the Jacobi explicit formula:
\[
[p/q]_{f,\z}(z)=\frac{\det\left|\begin{array}{cccc}
  (z-\z)^qS_{p-q}(z) & (z-\z)^{q-1}S_{p-q+1}(z) & \cdots & S_p(z) \\
  a_{p-q+1} & a_{p-q+2} & \cdots & a_{p+1} \\
  \vdots & \vdots &  & \vdots \\
  a_p & a_{p+1} & \cdots & a_{p+q} \\
\end{array}\right|}{\det\left|\begin{array}{cccc}
  (z-\z)^q & (z-\z)^{q-1} & \cdots & 1 \\
  a_{p-q+1} & a_{p-q+2} & \cdots & a_{p+1} \\
  \vdots & \vdots &  & \vdots \\
  a_p & a_{p+1} & \cdots & a_{p+q} \\
\end{array}\right|},
\]
with $S_k(z)=\left\{\begin{array}{cc}
  \dis\sum^k_{v=0}a_v(z-\z)^v, & k\ge0 \\
  0, & k<0. \\
\end{array}\right.$
For $\z=0$ we write $\cD_{p,q}$ instead of $\cD_{p,q}(0)$ and $[p/q]_f$ instead of $[p/q]_{f,0}$.

Next we present a theorem without proof,
concerning the Hankel determinant ($\ast$) of rational functions (see \cite{2}).
\begin{thm}\label{thm2.1}
Let $f=\frac{P}{Q}$ be a rational function, where $P,Q$ are polynomials, without any common factors.
If $\text{deg}P=\la$, $\text{deg}Q=\mi$ and $f$ is holomorphic around $\z\in\C$, then the following hold:

i)$f\in\cD_{\la,\mi}(\z)$,
ii)$f\in\cD_{p,\mi}(\z)$ for every $p>\la$,
iii)$f\in\cD_{\la,q}(\z)$ for every $q>\mi$ and
iv) $f\not\in\cD_{p,q}(\z)$, when $p>\la$ and $q>\mi$.
\end{thm}

If $\OO$ is an open set in $\C$ we denote by $H(\OO)$ the space of holomorphic functions
in $\OO$ endowed with the topology of uniform convergence on compacta. Hence, $H(\OO)$ is a Fr\'{e}chet space.
Therefore, Baire's Category theorem is at our disposal.

A holomorphic function $f\in H(\OO)$ belongs to $A^\infty(\OO)$ iff every derivative
$f^{(l)}$, $l=0,1,2,\ld$, extends continuously on $\oO$. We say that a sequence in $A^\infty(\OO)$,
$(f_m)_{m\in\N}$ converges $f_m\overset{A^\infty(\OO)}{\lra}f$, if and only if
$f^{(l)}_m\ra f^{(l)}$ uniformly on each compact subset of $\oO$, $\forall l$. Thus, $A^\infty(\OO)$ is
also a Fr\'{e}chet space and so Baire's theorem can be applied to any closed subspace of $A^\infty(\OO)$.
\section{Simply connected open sets}\label{sec3}
\noindent

Let $\OO\subseteq\C$ be a simply connected open set and $K,L$ two compact subsets of $\OO$.
Also, let $\cf\subseteq\N\times\N$ that contains a sequence $(\tp_n,\tq_n)$,
$n=1,2,\ld$ with $\tp_n\ra+\infty$. We define the following subsets of $H(\OO)$:
\begin{itemize}
\item $\cu(\OO,\cf,K,L)=\{f\in H(\OO)$: there exists a sequence $(p_n,q_n)_{n\in\N}$ in
$\cf$ such that $f\in \cD_{p_n,q_n}(\z)$, for all $\z\in L, n\in\N$ and
$\dis\sup_{\z\in L}\dis\sup_{z\in K}\big|[p_n/q_n]_{f,\z}(z)-f(z)\big|\lra0$, as $n\ra+\infty\}$.
\item $E(\OO,K,L,s,(p,q))=\{g\in H(\OO):g\in \cD_{p,q}(\z)$ for all $\z\in L$ and
$\dis\sup_{\z\in L}\dis\sup_{z\in K}\big|[p/q]_{g,\z}(z)-g(z)\big|<1/s\}$, where $s=1,2,\ld$ and $(p,q)\in \cf$.
\end{itemize}
\begin{rem}\label{rem3.1}
If $g\in E(\OO,K,L,s,(p,q))$ it follows that $[p/q]_{g,\z}$, $\z\in L$, does not have any poles in $K$.
The same holds for $[p_n/q_n]_{f,\z}$, $n\ge n_0$ (for some $n_0\in\N$), when $f\in \cu(\OO,\cf,K,L)$.
\end{rem}
\begin{prop}\label{prop3.2}
$\cu(\OO,\cf,K,L)=\dis\bigcap^\infty_{s=1}\dis\bigcup_{(p,q)\in \cf}E(\OO,K,L,s,(p,q))$.
\end{prop}
\begin{Proof}
It is standard and is omitted (see a similar proof in \cite{14}).
\qb
\end{Proof}
\begin{prop}\label{prop3.3}
The set $E(\OO,K,L,s,(p,q))$ is open in $H(\OO)$.
\end{prop}
\begin{Proof}
Let $g\in E(\OO,K,L,s,(p,q))$.
We consider another compact set $T\subseteq\OO$ such that $K\cup L\subseteq T^0$.
Hence, there exists $r>0$ such that for every $\z\in L$
it holds $\overline{\cD(\z,r)}\subseteq T^0$. Observe that the Hankel determinant $(\ast)$
and the coefficients of the numerator and the denominator of $[p/q]_{f,\z}$, depend polynomially
on $\frac{g^{(v)}(\z)}{v!}$, $v=0,1,\ld,p+q$ and continuously with respect to $\z\in L$.
This implies that for every $\e>0$ there exists $\tde>0$
such that for every $\tg\in H(\OO)$ with $\big|\frac{\tg^{(v)}(\z)}{v!}-\frac{g^{(v)}(\z)}{v!}\big|<\tde$,
$v=0,1,\ld,p+q$, $\z\in L$, it holds $\tg\in \cD_{p,q}(\z)$ and $\|[p/q]_{\tg,\z}-[p/q]_{g,\z}\|_{K}<\e$, $\forall \z\in L$.

If $0<\de<\min\{r^v\cdot\tde\;|\;v=0,1,\ld,p+q\}$ and
$\|\tg-g\|_T<\de$, then by Cauchy's estimates we obtain:
\[
\bigg|\frac{\tg^{(v)}(\z)}{v!}-\frac{g^{(v)}(\z)}{v!}\bigg|=
\bigg|\frac{(\tg-g)^{(v)}(\z)}{v!}\bigg|\le\frac{\|\tg-g\|_{\overline{D(\z,r)}}}{r^v}\le
\frac{\|\tg-g\|_T}{r^v}<\frac{\de}{r^v}<\tde.
\]
By choosing $\e$ and $\de$ properly we can achieve
$\dis\sup_{\z\in L}\dis\sup_{z\in K}\big|[p/q]_{\tg,\z}(z)-\tg(z)\big|<1/s$ and thus,
$\tg\in E(\OO,K,L,s,(p,q))$, which completes the proof.
\qb
\end{Proof}
\begin{prop}\label{prop3.4}
The set $\dis\bigcup_{(p,q)\in \cf}E(\OO,K,L,s,(p,q))$ is open and dense in $H(\OO)$, for all $s\in\N$.
\end{prop}
\begin{Proof}
Proposition \ref{prop3.3} implies that $\dis\bigcup_{(p,q)\in \cf}E(\OO,K,L,s,(p,q))$ is open.
According to Runge's theorem, the set of polynomials is
dense in $H(\OO)$, because $\OO$ is a simply connected open set. In order to show that
$\dis\bigcup_{(p,q)\in \cf}E(\OO,K,L,s,(p,q))$ is dense, let $P$ be a polynomial, $\e>0$ and $T\subseteq\OO$
a compact set. There exists $(p,q)\in \cf$ such that $p>\text{deg}P$.
We define $g(z)=P(z)+dz^p$, $d\in\C\sm\{0\}$.
It holds $\|g(z)-P(z)\|_T=\|dz^p\|_T<\e$, when $0<|d|<\frac{\e}{\|z\|^p_T}$.
Observe that for any $\z\in\C$ the development of $g$ is a polynomial of degree $p$ with $d$ as leading coefficient.
Therefore, one can easily see that, if $q\ge1$, the corresponding Hankel determinant
$(\ast)$ is equal to $(-1)^{\frac{q(q-1)}{2}}d^q\neq0$.
Thus, $g\in \cD_{p,q}(\z)$ for all $\z\in L$.
Since the Pad\'{e} approximant of $g$ is uniquely defined, it follows
that $[p/q]_{g,\z}(z)\equiv g(z)$, $\z\in L$, and so
$\dis\sup_{\z\in L}\dis\sup_{z\in K}\big|[p/q]_{g,\z}(z)-g(z)\big|=0<1/s$. The latter is also true in case $q=0$.
This implies $g\in\dis\bigcup_{(p,q)\in \cf}E(\OO,K,L,s,(p,q))$ and the proof is complete.
\qb
\end{Proof}
\begin{prop}\label{prop3.5}
Under the above assumptions and notation the set $\cu(\OO,\cf,K,L)$ is $G_\de$ and dense in $H(\OO)$.
[Hence, $\cu(\OO,\cf,K,L)\neq\emptyset$.]
\end{prop}
\begin{Proof}
It suffices to combine Proposition \ref{prop3.2} and Proposition \ref{prop3.4}
along with Baire's Category theorem.
\qb
\end{Proof}

Let us consider now $L$ to be a singleton $L=\{\z_0\}$. If we vary $K$ in an exhausting family $K_m$, $m\in\N$,
of compact subsets of $\OO$, then Baire's theorem implies that $\dis\bigcap^\infty_{m=1}\cu(\OO,\cf,K_m,\{\z_0\})$
is $G_\de$ and dense in $H(\OO)$. Hence, we get the following generic result.
\begin{thm}\label{thm3.6}
Under the above assumptions and notation
there exists a function $f\in H(\OO)$ with the following property:
"there exists a sequence $(p_n,q_n)\in\cf$, $n=1,2\ld$, such that
$f\in \cD_{p_n,q_n}(\z_0)$ for all $n$ and for every compact set $K\subseteq\OO$
it holds $\dis\sup_{z\in K}\big|[p_n/q_n]_{f,\z_0}(z)-f(z)\big|\lra0$, as $n\ra+\infty$". Furthermore, the set
of such functions is a $G_\de$ and dense subset of $H(\OO)$.
\end{thm}

Next we consider $K=L$ and vary it in an exhausting sequence $K_m$, $m\in\N$,
of compact subsets of $\OO$. By combining Proposition \ref{prop3.5} with Baire's theorem
we obtain the following generic property of holomorphic functions $f\in H(\OO)$.
\begin{thm}\label{thm3.7}
Let $\OO$ be a simply connected open set and $\cf\subseteq\N\times\N$ containing a sequence
$(\tp_n,\tq_n)\in\cf$, $n=1,2,\ld$, with $\tp_n\ra+\infty$. There exists a holomorphic function
$f\in H(\OO)$ with the following property:
"there exists a sequence $(p_n,q_n)\in\cf$, $n=1,2,\ld$, such that for every compact set $K\subseteq\OO$,
there is an $n(K)\in\N$ so that $f\in \cD_{p_n,q_n}(\z)$ for all $n\ge n(K)$, $\z\in L$ and
$\dis\sup_{\z\in K}\dis\sup_{z\in K}\big|[p_n/q_n]_{f,\z}(z)-f(z)\big|\lra0$, as $n\ra+\infty$". The set of such
functions $f\in H(\OO)$ is a $G_\de$ and dense subset of $H(\OO)$
endowed with the topology of uniform convergence on compacta.
\end{thm}
\begin{rem}\label{rem3.8}
Theorems \ref{thm3.6}, \ref{thm3.7} strengthen and simplify the proof of \textbf{Theorem 1} in \cite{3} and
also of a corollary of \textbf{Theorem 5.1} in \cite{5}.
\end{rem}
\section{Smooth functions on some simply connected open sets}\label{sec4}
\noindent

Let $\OO\subseteq\C$ be an open set such that $(\C\cup\infty)\sm\oO$ is connected; see also \cite{11}. We denote
$H(\oO)$ the set of holomorphic functions on some varying open set containing $\oO$. Let $X^\infty$ denote
the closure of $H(\oO)$ in $A^\infty(\OO)$. We do not know whether $X^\infty=A^\infty(\OO)$. A sufficient condition
in order for $X^\infty$ and $A^\infty(\OO)$ to coincide can be found in \cite{13}.
The condition states that there exists $M<\infty$, so that any two points in $\OO$
can be joined in $\OO$ by a curve $\Ga$ of length $|\Ga|\le M$.
\begin{lem}\label{lem4.1}
The set of polynomials is dense in $X^\infty$.
\end{lem}
\begin{Proof}
Let $f\in H(\oO)$, $T$ a compact subset of $\oO$, $\e>0$ and $N\in\N$. There exists $U\subseteq\C$ open such that $f\in H(U)$ and $\oO\subseteq U$. Also, there is $M>0$ such that $T\subseteq\oO\cap\overline{D(0,M)}$.
Observe that $(\C\cup\infty)\sm(\oO\cap\overline{D(0,M)})$ is connected.

\begin{itemize}
\item It is true that we can find $V\subseteq\C$ open, such that $\oO\cap\overline{D(0,M)}\subseteq V\subseteq U$
and $(\C\cup\infty)\sm V$ connected (equivalently $V$ is a simply connected open set) (\cite{8}, \cite{4}).
\end{itemize}

By Runge's theorem there exists a sequence of polynomials $(P_m)_{m\in\N}$,
such that $P_m\ra f$ uniformly on each compact subset of $V$. Weierstrass theorem implies that
$P^{(\ell)}_m\ra f^{(\ell)}$ uniformly on $T\subseteq\ V$, for every $\ell$. Hence, there exists $m_0\in\N$ such that
$\dis\|P^{(\ell)}_{m_0}-f^{(\ell)}\|_T<\e$, for all $0\le \ell\le N$. This completes the proof.
\qb
\end{Proof}

Let $\cf\subseteq\N\times\N$ be a subset containing a sequence $(\tp_n,\tq_n)$, $n\in\N$ with $\tp_n\ra+\infty$.
Also, let $K$ and $L$ two compact subsets of $\oO$. We define the following subsets of $X^\infty$:
\begin{itemize}
\item $\cb(\OO,\cf,K,L)=\{f\in X^\infty$: there exists a sequence $(p_n,q_n)\in\cf$, $n=1,2,\ld$ such that
$f\in\cD_{p_n,q_n}(\z)$ for all $\z\in L$, $n\in\N$ and
$\dis\sup_{\z\in L}\dis\sup_{z\in K}\big|[p_n/q_n]^{(\ell)}_{f,\z}(z)-f^{(\ell)}(z)\big|\lra0$,
as $n\ra+\infty$, for every $\ell=0,1,\ld$ \}.
\item $E^\infty(\OO,K,L,s,(p,q))=\{g\in X^\infty$: $g\in\cD_{p,q}(\z)$ for all $\z\in L$ and
$\dis\sup_{\z\in L}\dis\sup_{z\in K}\big|[p/q]^{(\ell)}_{g,\z}(z)-g^{(\ell)}(z)\big|<1/s$, for every $\ell=0,1,\ld,s$\},
where $s=1,2,\ld$ and $(p,q)\in\cf$.
\end{itemize}

We notice that when $\z\in bd(\OO)$, then $[p/q]_{f,\z}$ is the Pad\'{e} approximant of the series
$\dis\sum^\infty_{v=0}\frac{f^{(v)}(\z)}{v!}(z-\z)^v$, where $f^{(v)}(\z)$
is defined by continuity of the derivatives of $f$ in $\oO$.
\begin{prop}\label{prop4.2}
$\cb(\OO,\cf,K,L)=\dis\bigcap^\infty_{s=1}\dis\bigcup_{(p,q)\in \cf}E^\infty(\OO,K,L,s,(p,q))$.
\end{prop}
\begin{Proof}
It is standard and is omitted (see \cite{14}). \qb
\end{Proof}
\begin{prop}\label{prop4.3}
$E^\infty(\OO,K,L,s,(p,q))$ is open in $X^\infty$.
\end{prop}
\begin{Proof}
Let $g\in E^\infty(\OO,K,L,s,(p,q))$. Observe that the Hankel determinant ($\ast$) and the coefficients
of the numerator and the denominator of $[p/q]^{(\ell)}_{f,\z}$, depend polynomially on $\frac{g^{(v)}(\z)}{v!}$,
$v=0,1,\ld,p+q$, for every $\ell=0,1,\ld$, and continuously with respect to $\z\in L$.
This implies that for every $\e>0$ there exists $\tde>0$ such that
for every function $\tg\in X^\infty$ with $\big|\frac{\tg^{(v)}(\z)}{v!}-\frac{g^{(v)}(\z)}{v!}\big|<\tde$,
$v=0,1,\ld,p+q$, $\z\in L$, it holds $\tg\in\cD_{p,q}(\z)$ and
$\big\|[p/q]^{(\ell)}_{\tg,\z}(z)-[p/q]^{(\ell)}_{g,\z}(z)\big\|_K<\e$, for all $\z\in L$ and $\ell=0,1,\ld,s$.

We consider T a compact subset of $\oO$ such that $K\cup L\subseteq T$.
If $0<\de<\min\{r^v\cdot\tde\;|\;v=0,1,\ld,p+q\}$,
then for every function $\tg\in X^\infty$ satisfying $\|\tg^{(\ell)}-g^{(\ell)}\|_T<\de$,
for all $0\le \ell\le\max\{p+q,s\}$, it holds
\[
\bigg|\frac{\tg^{(v)}(\z)}{v!}-\frac{g^{(v)}(\z)}{v!}\bigg|=
\frac{\big|\tg^{(v)}(\z)-g^{(v)}(\z)\big|}{v!}\le\frac{\|\tg^{(v)}-g^{(v)}\|_L}{v!}\le
\frac{\|\tg^{(v)}-g^{(v)}\|_T}{v!}<\frac{\de}{v!}<\tde.
\]
By choosing $\e$ and $\de$ properly we obtain
$\dis\sup_{\z\in L}\dis\sup_{z\in K}\big|[p/q]^{(\ell)}_{\tg,\z}(z)-\tg^{(\ell)}(z)\big|<1/s$, for every $\ell=0,1,\ld,s$
and thus, $\tg\in E^\infty(\OO,K,L,s,(p,q))$, which completes the proof. \qb
\end{Proof}
\begin{prop}\label{prop4.4}
The set $\dis\bigcup_{(p,q)\in\cf}E^\infty(\OO,K,L,s,(p,q))$ is open and dense in $X^\infty$.
\end{prop}
\begin{Proof}
By Proposition \ref{prop4.3} this union is open in $X^\infty$. According to Lemma \ref{lem4.1},
polynomials are dense in $X^\infty$. Let $P$ be a polynomial, $T$ a compact subset of $\oO$,
$N\in\N$ and $\e>0$. There exists $(p,q)\in\cf$ such that $p>\text{deg}P$. We define $g(z)=P(z)+dz^p$, $d\in\C\sm\{0\}$.
It is clear that $g\in X^\infty$ and
$\|g^{(\ell)}(z)-P^{(\ell)}(z)\|_T=|d|\cdot\|(z^p)^{(\ell)}\|_T<\e$, for all $\ell=0,1,\ld,N$,
when $0<|d|<\frac{\e}{\max_{0\le \ell\le N}\|(z^p)^{(\ell)}\|_T}$.
By developing $g$ around $\z\in\C$, we can easily see that, if $q\ge1$,
the corresponding Hankel determinant $(\ast)$ is independent of $\z$ and equal to $(-1)^{\frac{q(q-1)}{2}}\cdot d^q\neq0$.
Thus, $g\in\cD_{p,q}(\z)$ and $[p/q]_{g,\z}(z)\equiv g(z)$ for all $\z\in L$. The latter also holds in case $q=0$.
Hence, it follows that 
$\dis\sup_{\z\in L}\dis\sup_{z\in K}\big|[p/q]^{(\ell)}_{g,\z}(z)-g^{(\ell)}(z)\big|=0<1/s$, for every $l=0,\ld,s$ 
and $g\in\dis\bigcup_{(p,q)\in\cf}E^\infty(\OO,K,L,s,(p,q))$, which completes the proof.
\qb
\end{Proof}
\begin{prop}\label{prop4.5}
Under the above assumptions and notation $\cb(\OO,\cf,K,L)$ is $G_\de$ and dense in $X^\infty$.
[Hence, $\cb(\OO,\cf,K,L)\neq\emptyset$.]
\end{prop}
\begin{Proof}
The proof follows immediately, if we combine Proposition \ref{prop4.2} and Proposition \ref{prop4.4} along with
Baire's Category theorem.
\qb
\end{Proof}

Let $L=\{\z_0\}$, $\z_0\in\oO$ and let $K$ vary in the sequence $K_m=\oO\cap\overline{D(o,m)}$, $m\in\N$.
Then Proposition \ref{prop4.5} with Baire's theorem implies that
$\dis\bigcap^\infty_{m=1}\cb(\OO,\cf,K_m,\{\z_0\})$ is $G_\de$ and dense in $X^\infty$.
This yields the following generic result.
\begin{thm}\label{thm4.6}
Under the above assumptions and notation
there exists a function $f\in A^\infty(\OO)$ and a sequence $(p_n,q_n)\in\cf$, $n=1,2,\ld$, such that
for every compact set $K\subseteq\oO$ the following hold: "$f\in\cD_{p_n,q_n}(\z_0)$ for all $n=1,2,\ld$ and
$\dis\sup_{z\in K}\big|[p/q]^{(\ell)}_{f,\z_0}(z)-f^{(\ell)}(z)\big|\lra0$, as $n\ra+\infty$, for every $\ell$".
Moreover, the set of such functions is $G_\de$ and dense in $X^\infty$,
the closure of polynomials in $A^\infty(\OO)$.
\end{thm}

Next, we let $K=L$ and vary it in the sequence $K_m=\oO\cap\overline{D(0,m)}$, $m\in\N$.
Baire's theorem together with Proposition \ref{prop4.5} implies that
$\dis\bigcap^\infty_{m=1}\cb(\OO,\cf,K_m,K_m)$ is $G_\de$ and dense in $X^\infty$. The latter is
expressed in the following.
\begin{thm}
Generically every function $f\in X^\infty$ has the following property:
"there exists a sequence $(p_n,q_n)_{n\in\N}\in\cf$ such that for every compact set
$K\subseteq\oO$ there is an $n(K)\in\N$ so that for every $n\ge n(K)$ it holds $f\in\cD_{p_n,q_n}(\z)$ for all $\z\in K$ and
$\dis\sup_{\z\in K}\dis\sup_{z\in K}\big|[p/q]^{(\ell)}_{f,\z}(z)-f^{(\ell)}(z)\big|\lra0$,
as $n\ra+\infty$, for every $\ell$".
\end{thm}
\section{Arbitrary open sets}\label{sec5}
\noindent

Let $\OO\subseteq\C$ be an arbitrary open set and let $K,L$ be two compact subsets of $\OO$. Also, let
$\cf\subseteq\N\times\N$ which contains a sequence $(\tp_n,\tq_n)\in\cf$, $n\in\N$, such that $\tp_n\ra+\infty$
and $\tq_n\ra+\infty$, as $n\ra+\infty$. We define the following subsets of $H(\OO)$:
\begin{itemize}
\item $\cu(\OO,\cf,K,L)=\{f\in\ H(\OO)$: there exists a sequence $(p_n,q_n)_{n\in\N}$ in $\cf$ such that
$f\in\cD_{p_n,q_n}(\z)$ for all $n\in\N$, $\z\in L$ and
$\dis\sup_{\z\in L}\dis\sup_{z\in K}\big|[p_n/q_n]_{f,\z}(z)-f(z)\big|\lra0$, as $n\ra+\infty\}$.
\item $E(\OO,K,L,s,(p,q))=\{g\in H(\OO):g\in \cD_{p,q}(\z)$ for all $\z\in L$ and
$\dis\sup_{\z\in L}\dis\sup_{z\in K}\big|[p/q]_{g,\z}(z)-g(z)\big|<1/s\}$, where $s=1,2,\ld$ and $(p,q)\in \cf$.
\end{itemize}
\begin{prop}\label{prop5.1}
$\cu(\OO,\cf,K,L)=\dis\bigcap^\infty_{s=1}\dis\bigcup_{(p,q)\in\cf}E(\OO,K,L,s,(p,q))$.
\end{prop}
\begin{prop}\label{prop5.2}
The set $E(\OO,K,L,s,(p,q))$ is open in $H(\OO)$.
\end{prop}
\begin{Proof}
The proof is the same with that of Proposition \ref{prop3.3} and is omitted.\qb
\end{Proof}
\begin{prop}\label{prop5.3}
The set $\dis\bigcup_{(p,q)\in\cf}E(\OO,K,L,s,(p,q))$ is open and dense in $H(\OO)$, for all $s\in\N$.
\end{prop}
\begin{Proof}
Proposition \ref{prop5.2} implies that this union is open. According to Runge's theorem, the rational
functions with poles off $\OO$ are dense in $H(\OO)$. In order to show that
$\dis\bigcup_{(p,q)\in\cf}E(\OO,K,L,s,(p,q))$ is dense, let $R=\frac{A}{B}$ be a rational function without
poles in $\OO$, where $A,B$ are polynomials without any common factors.
Also, let $\e>0$ and $T\subseteq\OO$ a compact set.
\begin{itemize}
  \item There exists $(p,q)\in F$ such that $p>\text{deg}A$ and $q>\text{deg}B$.
\end{itemize}
We define $g(z)=\frac{A(z)+dz^p}{B(z)}$, where $d\in\C\sm\{0\}$.
It holds $\dis\inf_{z\in T}|B(z)|>0$, because $R$ has no poles in $\OO$. Hence, we have
$\|g(z)-R(z)\|_T=\big\|\frac{dz^p}{B(z)}\big\|_T<\e$, when $0<|d|<\frac{\e\cdot\inf_T|B(z)|}{\|z\|^p_T}$.

Furthermore, we can choose $d\in\C\sm\{0\}$ satisfying the above so that $g$ is irreducible.
This can be achieved if $d\neq-\frac{A(\rho_i)}{\rho^p_i}$, $i=1,\ld,\ti$, where $\rho_i$ are
the zeros of $B$. We notice that if $\rho_{i_0}=0$, $i_0\in\{1,\ld,\ti\}$, then
$A(\rho{i_0})\neq0$ or $A$ and $B$ would have had a common factor.
Thus, according to Theorem \ref{thm2.1} we have $g\in\cD_{p,q}(\z)$, for all $\z\in L$,
since $p=\text{deg}(A(z)+dz^p)$ and $q>\text{deg}B(z)$. It follows that $[p/q]_{g,\z}\equiv g$,
for every $\z\in L$ and so $\dis\sup_{\z\in L}\dis\sup_{z\in K}\big|[p/q]_{g,\z}(z)-g(z)\big|=0<1/s$.
This implies that $g\in\dis\bigcup_{(p,q)\in\cf}E(\OO,K,L,s,(p,q))$ and the proof is complete.
\qb
\end{Proof}
\begin{prop}\label{prop5.4}
Under the above assumptions and notation the set $\cu(\OO,\cf,K,L)$ is
$G_\de$ and dense in $H(\OO)$. [Hence, $\cu(\OO,\cf,K,L)\neq\emptyset$.]
\end{prop}
\begin{Proof}
It suffices to combine Proposition \ref{prop5.1} and Proposition \ref{prop5.3} along with
Baire's Category theorem.
\qb
\end{Proof}

Special choices of $K$ and $L$ combined with Proposition \ref{prop5.4} and Baire's theorem
yield the following generic results.
\begin{thm}\label{thm5.5}
Under the above assumptions and notation, let $\z_0\in\OO$. Then there exists $f\in H(\OO)$
with the following property: "there exists a sequence $(p_n,q_n)_{n\in\N}$ in $\cf$ such that
$f\in\cD_{p_n,q_n}(\z_0)$, $\forall n\in\N$ and for every  compact set $K\subseteq\OO$ it holds
$\dis\sup_{z\in K}\big|[p_n/q_n]_{f,\z_0}(z)-f(z)\big|\lra0$, as $n\ra+\infty$". In addition,
the set of such functions is a $G_\de$ and dense subset of $H(\OO)$.
\end{thm}
\begin{thm}
Let $\OO\subseteq\C$ be an arbitrary open set and $\cf\subseteq\N\times\N$ containing a
sequence $(\tp_n,\tq_n)\in\cf$, $n=1,2,\ld$ with $\tp_n\ra+\infty$ and $\tq_n\ra+\infty$, as $n\ra+\infty$.
Then there exists a holomorphic function $f\in H(\OO)$ so that the following hold:
"there exists a sequence $(p_n,q_n)_{n\in\N}$ in $\cf$ such that for every $K\subseteq\OO$ compact, there is
an $n(K)\in\N$ so that $f\in\cD_{p_n,q_n}(\z)$ for all $n\ge n(K)$, $\z\in K$ and
$\dis\sup_{\z\in K}\dis\sup_{z\in K}\big|[p_n/q_n]_{f,\z}(z)-f(z)\big|\lra0$, as $n\ra+\infty$". Further,
the set of such functions $f\in H(\OO)$ is $G_\de$ and dense in $H(\OO)$
endowed with the topology of uniform convergence on compacta.
\end{thm}
\section{Smooth functions on arbitrary open sets}\label{sec6}
\noindent

Let $\OO\subseteq\C$ be an open set. We denote by $H(\oO)$ the set of holomorphic functions
on some varying open set containing $\oO$. Let $X^\infty$ denote the closure of $H(\oO)$ in
$A^\infty(\OO)$. By Runge's theorem the rational functions with poles off $\oO$ are dense in $X^\infty$.

Let $\cf\subseteq\N\times\N$ which contains a sequence
$(\tp_n,\tq_n)\in\cf$, $n=1,2,\ld$ with $\tp_n\ra+\infty$ and $\tq_n\ra+\infty$, as $n\ra+\infty$. Also,
let $K$ and $L$ be two compact subsets of $\oO$. We define the following subsets of $X^\infty$:
\begin{itemize}
\item $\cb(\OO,\cf,K,L)=\{f\in X^\infty$: there exists a sequence $(p_n,q_n)\in\cf$, $n=1,2,\ld$ such that
$f\in\cD_{p_n,q_n}(\z)$ for all $\z\in L$, $n\in\N$ and
$\dis\sup_{\z\in L}\dis\sup_{z\in K}\big|[p_n/q_n]^{(\ell)}_{f,\z}(z)-f^{(\ell)}(z)\big|\lra0$,
as $n\ra+\infty$, for every $\ell=0,1,\ld$ \}.
\item $E^\infty(\OO,K,L,s,(p,q))=\{g\in X^\infty$: $g\in\cD_{p,q}(\z)$ for all $\z\in L$ and
$\dis\sup_{\z\in L}\dis\sup_{z\in K}\big|[p/q]^{(\ell)}_{g,\z}(z)-g^{(\ell)}(z)\big|<1/s$, for every $\ell=0,1,\ld,s$\},
where $s=1,2,\ld$ and $(p,q)\in\cf$.
\end{itemize}

We notice that when $\z\in bd(\OO)$, then $[p/q]_{f,\z}$ is the Pad\'{e} approximant of the series
$\dis\sum^\infty_{v=0}\frac{f^{(v)}(\z)}{v!}(z-\z)^v$, where $f^{(v)}(\z)$ is defined by continuity.
\begin{prop}\label{prop6.1}
$\cb(\OO,\cf,K,L)=\dis\bigcap^\infty_{s=1}\dis\bigcup_{(p,q)\in \cf}E^\infty(\OO,K,L,s,(p,q))$.
\end{prop}
\begin{prop}\label{prop6.2}
$E^\infty(\OO,K,L,s,(p,q))$ is open in $X^\infty$.
\end{prop}
\begin{Proof}
It is similar to the proof of Proposition \ref{prop4.3}. \qb
\end{Proof}
\begin{prop}\label{prop6.3}
The set $\dis\bigcup_{(p,q)\in\cf}E^\infty(\OO,K,L,s,(p,q))$ is open and dense in $X^\infty$.
\end{prop}
\begin{Proof}
By Proposition \ref{prop6.2} this union is open in $X^\infty$. Let $R=\frac{A}{B}$ be a rational
function with poles off $\oO$, where $A,B$ are polynomials without any common factors. Hence,
there exists an open set $U\supseteq\oO$ such that $R\in H(U)$. Also, let
$T$ be a compact subset of $\oO$, $N\in\N$ and $\e>0$.
\begin{itemize}
\item There exists $(p,q)\in F$ such that $p>\text{deg}A$ and $q>\text{deg}B$.
\end{itemize}
We define $g_m(z)=\frac{A(z)+d_mz^p}{B(z)}\in X^\infty$, where $d_m\in\C\sm\{0\}$, $m\in\N$.
If $\rho_1,\ld,\rho_\ti$ are the zeros of $B$, then we consider $d_m\neq\-\frac{A(\rho_i)}{\rho^p_i}$,
for all $i=1,\ld,\ti$, $\forall m\in\N$.

Thus, if we let $d_m\ra0$, as $m\ra+\infty$, then $g_m\ra R$ uniformly on each compact subset
of $U$. By Weierstrass' theorem we have $g^{(\ell)}_m\ra R^{(\ell)}$ uniformly on $T$, for every $\ell$. It follows
that there exists $m_0\in\N$ such that $\|g^{(\ell)}_{m_0}(z)-R^{(\ell)}(z)\|_T<\e$, for all $\ell=0,\ld,N$.

Observe that $g_{m_0}$ is irreducible. According to Theorem \ref{thm2.1}, it holds
$g_{m_0}\in\cD_{p,q}(\z)$, for all $\z\in L$,
since $p=\text{deg}(A(z)+d_{m_0}z^p)$ and $q>\text{deg}B(z)$. It follows that
$[p/q]_{g_{m_0},\z}\equiv g_{m_0}$, for all $\z\in L$ and
$\dis\sup_{\z\in L}\dis\sup_{z\in K}\big|[p/q]^{(\ell)}_{g_{m_0},\z}(z)-g^{(\ell)}_{m_0}(z)\big|=0<1/s$,
for all $\ell=0,1,\ld,s$. Hence, $g\in\dis\bigcup_{(p,q)\in\cf}E^\infty(\OO,K,L,s,(p,q))$ which
completes the proof. \qb
\end{Proof}
\begin{prop}\label{prop6.4}
Under the above assumptions and notation $\cb(\OO,\cf,K,L)$ is $G_\de$ and dense in $X^\infty$.
[Hence, $\cb(\OO,\cf,K,L)\neq\emptyset$.]
\end{prop}
\begin{Proof}
The proof follows by combining Propositions \ref{prop6.1}, \ref{prop6.3} with Baire's
Category theorem. \qb
\end{Proof}

By varying $K$,$L$ and applying Baire's theorem, we obtain the following.

\begin{thm}\label{thm6.5}
Under the above assumptions and notation, let $\z_0\in\oO$. Then there exists a
function $f\in A^\infty(\OO)$ and a sequence $(p_n,q_n)\in\cf$, $n=1,2,\ld$, such that for
every compact set $K\subseteq\oO$ the following hold: "$f\in\cD_{p_n,q_n}(\z_0)$ for all $n=1,2\ld$
and $\dis\sup_{z\in K}\big|[p/q]^{(\ell)}_{f,\z_0}(z)-f^{(\ell)}(z)\big|\lra0$, as $n\ra+\infty$, for every $\ell$".
Moreover, the set of such functions is $G_\de$ and dense in $X^\infty$.
\end{thm}
\begin{thm}\label{thm6.6}
Let $\OO\subseteq\C$ be an open set. Also, let $\cf\subseteq\N\times\N$ containing a sequence
$(\tp_n,\tq_n)_{n\in\N}$ with $\tp_n\ra+\infty$ and $\tq_n\ra+\infty$, as $n\ra+\infty$.
There exists a function $f\in A^\infty(\OO)$ with the following property:
"there exists a sequence $(p_n,q_n)_{n\in\N}\subseteq\cf$, such that for every compact set
$K$, there is an $n(K)\in\N$ so that for all $n\ge n(K)$ it holds $f\in\cD_{p_n,q_n}(\z)$,
for every $\z\in K$ and
$\dis\sup_{\z\in K}\dis\sup_{z\in K}\big|[p_n/q_n]^{(\ell)}_{f,\z}(z)-f^{(\ell)}(z)\big|\lra0$, as $n\ra+\infty$,
for every $\ell$". The set of such functions is $G_\de$ and dense in $X^\infty$.
\end{thm}

\end{document}